\documentclass{amsart}

\usepackage{amssymb}
\usepackage{graphicx}

\usepackage[margin=1in]{geometry}  % set the margins to 1in on all sides
\usepackage{graphicx}              % to include figures
\usepackage{amsmath}               % great math stuff
\usepackage{amsfonts}              % for blackboard bold, etc
\usepackage{amsthm}                % better theorem environments

\begin{document}

\nocite{*}
\title{An analysis of IQ LINK$^{\textregistered}$ }

\author{Donna A. Dietz\\
Department of Mathematics and Statistics\\
 American University\\ Washington, DC, USA}
\email{dietz@american.edu}

\maketitle

\begin{abstract}
This is a theoretical and computational strategy exploration of the
visually attractive game IQ-Link.  Not all games which are visually
appealing are worthy of your time as a puzzlist. This analysis gives a
would-be addict some idea of what type of game they are falling for.
\end{abstract}

\section {Introduction}
From the moment I first opened my new IQ-Link puzzle (created by Raf
Peeters), I felt a very strong attraction to the toy.  As Raf says,
``It almost seems like the puzzle pieces are jewels or are made of
candy... The object of the game is to make all puzzle pieces fit on
the game board.''  The pieces themselves are colorful, shiny,
translucent, and smooth to the touch.  They even resonate with a
gently musical percussive sound on their board!  The toy is portable
and inexpensive, so it makes a great gift.

But what is it, really? Should you bother?  You get a sense that
you're trying to solve a problem that has connections to molecular
chemistry, but what if you're just being fooled by its superficial
beauty? Perhaps you, like me, feel the Sudoku ``strategy'' of
trial-and-error-until-you-drop is not really strategy at
all. According to computer security expert Ben Laurie, ``Sudoku is a
denial of service attack on human intellect''\cite{norvig}.  Peter Norvig, in an attempt to cure his wife from this
``virus'', wrote an article on his website, ``Solving Every Sudoku
Puzzle'' \cite{norvig}.  His code fits on a single page and can
solve the hardest Sudoku puzzle he could find in only 0.01 seconds.  What
if IQ-Link is a game of that sort?

In the puzzles of yesteryear, puzzles were created by hand and
intended to be solved by hand.  Now, computers can generate endless
``puzzles'' for humans, but are they really anything other than
clerical speed and accuracy drills? When we, as recreational
mathematicians, speak to this point, we can give puzzlists a general
idea about new games when they come out on the market.  This gives a
sense of who would like the toy. In this discussion, I shed
light on this complexity question for IQ-Link, a game for which I have
found no theoretical research.

Ultimately, the puzzles were all solved algorithmically after several
theoretical insights about the puzzle which led to algorithmic
adjustments. These insights are often the same insights that shed
light on how a puzzlist would experience the puzzle.  They brought the
complexity down to a computable level, but a few puzzles still took a
great deal of computation time (days).  It is my opinion that this toy
is well worth the time spent playing with it, as the puzzlist must
continue to develop improvements to their strategies at each
difficulty level.  These are true strategic improvements and not
merely the search through more combinations.

\section{Background}

IQ-Link is a collection of puzzles or one-person games, also called
challenges, played on a board with 12 pieces.
IQ-Link analysis lies in the study of combinatorial game theory.
There are no dice or
other random events or any hidden information in the challenges. Each
challenge begins with some of the pieces already being permanently
placed on the board, and the goal is to fit all the remaining pieces
onto the board without collisions.

In order to place this into its proper context, I will introduce some
terminology, but it is not necessary to read this background section
in order to understand the remainder of the talk/paper.
There are some variations in definitions for what consititutes a game,
and whether or not a one-player game is a game or not. I prefer less
restrictive definitions.  By the definitions set forth in Valia
Mitsou's PhD Thesis \cite{mitsou}, IQ-Link would be classified as a
one-player game, or a puzzle.  Similarly, there is debate about what
constitutes a combinatorial game,
but I'm going with the combinatorial game definition in {\it Winning
  Ways} which states that for a game to be combinatorial, both players
have to know what is going on.  ``There is complete information...
There is no occasion for bluffing'' \cite{winning}. 

Even though there
is no randomness in IQ-Link, variations could easily be
made to
give it randomness or even turn it into 
a two-player game.
I like the the definition given by Harish Karnick \cite{karnick}
that you can have a player called ``Nature'' which can introduce
random acts. This keeps the definition of combinatorial game as broad as possible.

This talk/paper is partly about theoretical insights, but it reflects mostly on
results from running algorithms on computers. Perhaps that
seems a little bit like cheating.  But is it?  To get a better idea of
what other researchers are doing, and thus what they feel is fair, I
want to go back to some ideas in the introduction of Valia Mitsou's
thesis \cite{mitsou}. She
states, ``the subject of this thesis is the algorithmic properties of
one-player games people enjoy playing'' and asks, ``can we design
efficient computer programs that... solve any instance of the puzzle
in question?''  She's not the only one asking such questions.

Interest in this area has been exploding recently.  According to the
editors of a January 2015 Dagstahl report \cite{dagstahl}, there is a
great need for interdisciplinary research alliances in order to
advance the boundaries of our current understanding of both combinatorial games, but also other types of games. Towards this goal,
they set up the conference, ``Artificial and computational
intelligence in games: integration'', with multiple workgroups.  One
such workgroup was ``Algorithms that learn to play like people,'' led
by Julian Togelius.  This interest by the conference organizers
indicates that, it is an active area of research to
create algorithms that emulate people. This is done  partly to classify the
difficulty of computer-generated puzzles (such as done by a team of
researchers including Mark van Kreveld \cite{kreveld} and
\cite{kreveld2}), but also sometimes to create automata to play
against humans.  Computers can simulate human dexterity (or reaction
time) as well as human strategy \cite{isaksen}.  IQ-Link is a game of
pure strategy like Chess, Go, or Suduko, and would appear in one corner of Aaron Isaksen's \cite{isaksen} dexterity / strategy
classification diagram. It would be in the opposite corner as ``Flappy Bird'' which is all about
dexterity and requires little to no strategy. As I've done here, many researchers focus
entirely on games of pure strategy. For example,
in a 2002 volume of studies edited by
Richard Nowakowski, Aviezre Fraenkel
\cite{fraenkel} and others do this.

In Robert Hearn's 1987 PhD thesis \cite{hearn} (under the guidance of Erik
Demaine and Gerald Sussman), he states that, ``There is a fundamental connection between the notions of game an computation... and
indeed... games should be thought of as a valid model of computation,
just as Turing machines are.''  He goes on to state that ``another
curious property of games is that they tend to be as hard as possible
given its general characteristics''.  This is reassuring to
me as it is  consistent with my experiences in attempting to
create and refine algorithms to solve IQ-Link.

In researching these various papers and talks, I see that many
researchers have attacked many games over the years.  In this
paper/talk, I share a partial understanding of
IQ-Link.
I have seen nothing published about this
particular game yet, so I feel it's a reasonable addition to an
ongoing series of discussions about  no-chance games and
puzzles in ``puzzleland''.

\section{Methods}
Just because it's possible to
write code to solve a puzzle does not mean it's easy to solve by hand.
Something which is somehow easy on a computer may be hard by hand and
possibly vice versa.  While constructing an algorithm, it naturally makes sense to keep playing with the puzzle
by hand, and to adjust the algorithm to imitate the by-hand attempts.
Of course this process will never fully mimic a human, but the process
has to keep feeding back on itself if we are to use the computer to
help predict other puzzlists' reactions to the puzzle.

One thing that seemed easy by hand but time consuming to compute was the
linking condition on the pieces which will be described in Section
\ref{linkage}.  This was not initially in my algorithm but was vital
to solving the more advanced puzzles.  The fact that a solution needs
to be unique is a big hint when solving by hand.  Again, the ultimate
goal for me was to find out if the higher levels of this puzzle will continue
to be fun for someone who doesn't like Sudoku-like games. This was done
by gradually improving the strategy and seeing how much they improved
the search algorithm for the harder puzzles.

\section{Toy mechanics}
\label{mechanics}
\begin{figure}
\begin{center}
\includegraphics[width=3in]{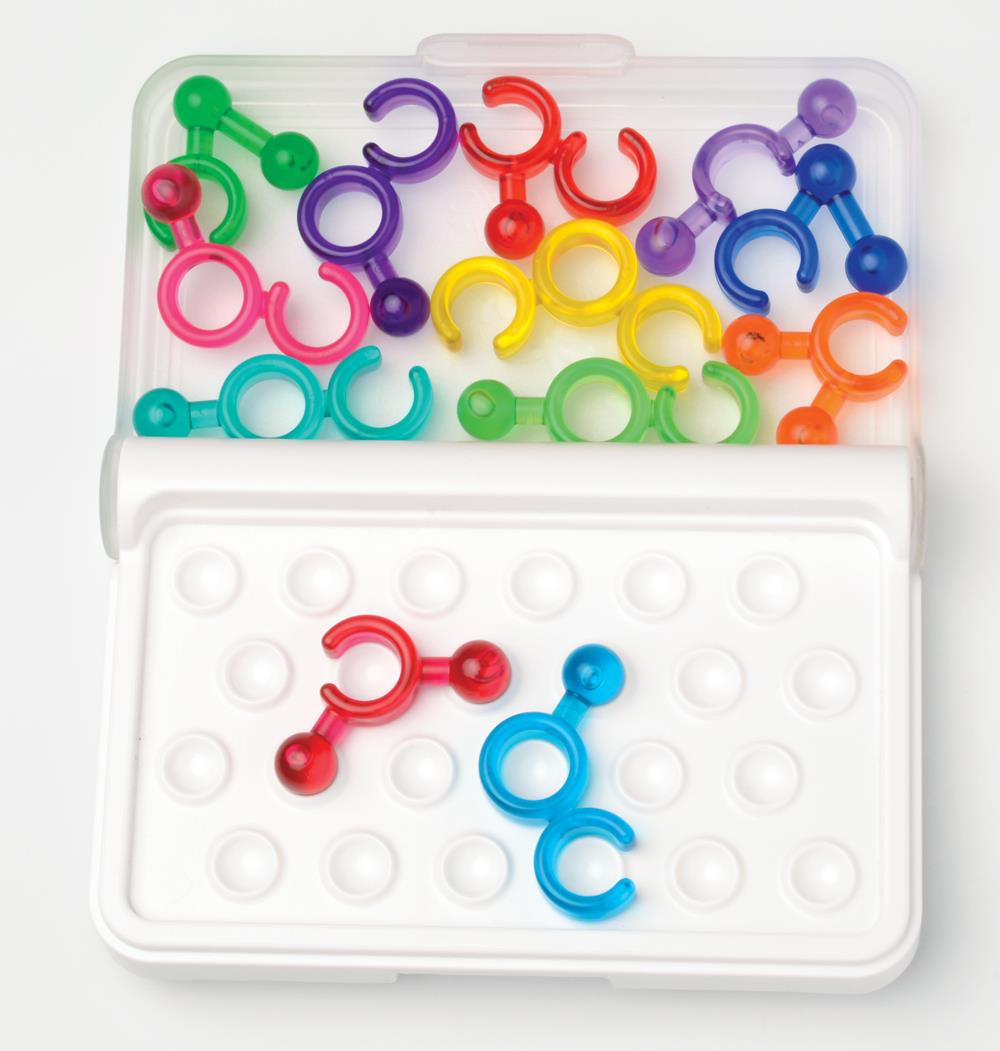}\\
\end{center}
\caption{This is the game board, box, and pieces.}
\label{gameboard}
\end{figure}

\begin{figure}
\begin{center}
\includegraphics[width=6in]{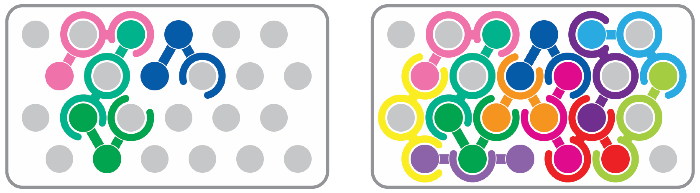}\\
\end{center}
\caption{This is Puzzle 92 and its solution.}
\label{pubpuz92}
\end{figure}

The game (whose manufacturer is found online \cite{smartgames})
is shown in Figure \ref{gameboard} as it appears when a puzzlist
begins to play with it. The puzzle comes with 120 challenges and
includes the solution key.  One challenage and its solution are shown
in Figure \ref{pubpuz92}. Pieces must rest in the raised spots on the
board in order for the placements to be legal. The entire puzzle must
be linked to form a connected graph, and the solutions to the given
puzzles are unique.

\section{Puzzle difficulty levels}
In the published puzzles (which come with an answer key), there are 24
puzzles each in 5 difficulty levels: Starter, Junior, Expert, Master,
and Wizard. This makes 120 puzzles.  See Table \ref{puzzlelevels}. All puzzles were solved by algorithms explained in this talk/paper and verified against the published key.

\begin{table}[]
  \begin{center}
  \caption{Puzzle difficulty levels}
  \label{puzzlelevels}
\begin{tabular}{|l|l|l|l|}
\hline
\textbf{puzzle level} & \textbf{puzzle numbers} & \textbf{number of pieces given} & \textbf{number of pieces to place} \\ \hline
Starter               & 1-24                    & 9-10                            & 2-3                                \\ \hline
Junior                & 25-48                   & 5-9                             & 3-7                                \\ \hline
Expert                & 49-72                   & 4-7                             & 5-8                                \\ \hline
Master                & 73-96                   & 3-6                             & 6-9                                \\ \hline
Wizard                & 97-120                  & 2-4                             & 8-10                               \\ \hline
\end{tabular}
\end{center}
\end{table}

\section{A hexagonal model}
\begin{figure}
\begin{center}
\includegraphics[width=3in]{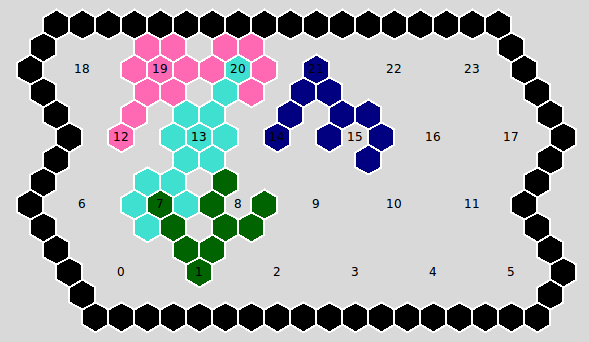}\\
\end{center}
\caption{A challenge.}
\label{92Puzzle}
\end{figure}

\begin{figure}
\begin{center}
\includegraphics[width=3in]{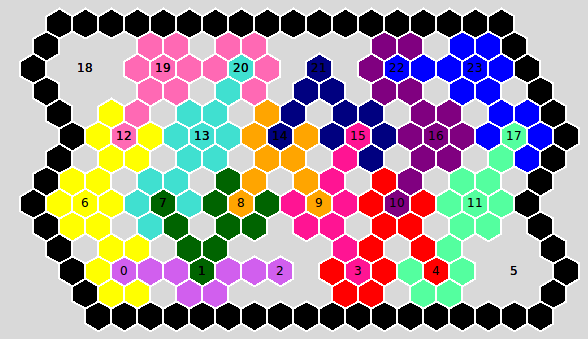}\\
\end{center}
\caption{The solution.}
\label{P92_solved}
\end{figure}

The interactions of the pieces and legal placements of the pieces on
the board suggest a hexagonal tiling model. Since this simplification
is easier for the computerized analysis, this will be used throughout
the article. The same puzzle and solution shown in Figure
\ref{pubpuz92} are again shown in Figures \ref{92Puzzle} and
\ref{P92_solved} but in hexagonal form.

\begin{figure}
\begin{center}
\includegraphics[width=5in]{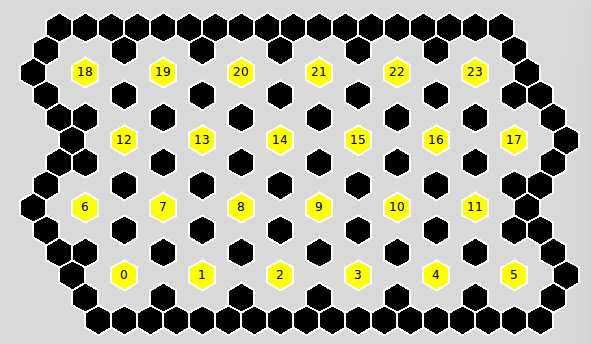}\\
\end{center}
\caption{The game board with anchor spots.}
\label{ReachablesAndHomes}
\end{figure}

The pieces can be grouped into
families by how many ways they may be legally placed on the board.
In the hexagonal diagram of the game board (Figure \ref{ReachablesAndHomes}), there are 24 special spots which act as
``anchors'' (colored yellow).  In the actual toy, these are physically raised circles
(called ``spaces'' in the instructions) which will support either a ball
on the inside, or a ring on the outside, or both- when the ring is not
a full ring.  In the instructions, we are told that not all ``spaces''
must be used in order to achieve a solution.  We are warned, however,
that the entire graph of the solution must be linked together.  The
ball-and-circle joints are actually vertices in a graph, and the
graph is connected in a legal solution.

\begin{table}[]
    \begin{center}
  \caption{Puzzle pieces}
  \label{puzzlepieces}
\begin{tabular}{|ll|llll|}
\hline
\textbf{placements} & \textbf{and family names}    & \textbf{diagrams} & &&  \\ \hline
72    & straight &   \includegraphics[width=1in]{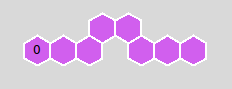}    & light purple  &&  \\ \hline
144         & straight   &   \includegraphics[width=1in]{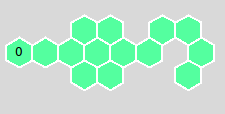}       & light green   &  \includegraphics[width=1in]{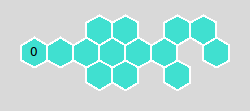}     &  aqua      \\ \hline
     &      &     \includegraphics[width=1in]{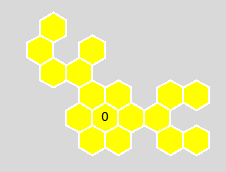}  & yellow   &     \includegraphics[width=1in]{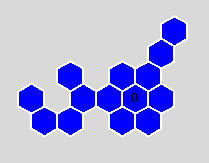}  & blue \\ 
    &      &     \includegraphics[width=1in]{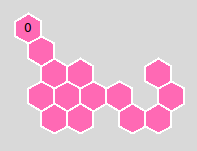}  & light pink   &     \includegraphics[width=1in]{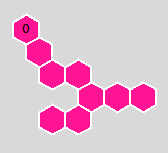}  & dark pink \\ 
152 & obtuse     &     \includegraphics[width=1in]{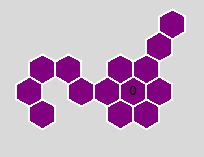}  & dark purple && \\ \hline
                          &        &  \includegraphics[width=1in]{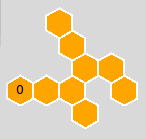}    & orange   &  \includegraphics[width=1in]{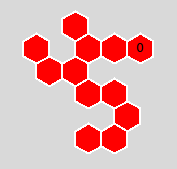}    & red   \\ 
180&  acute    &  \includegraphics[width=1in]{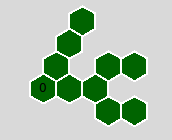}    & dark green  &  \includegraphics[width=1in]{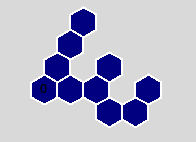}    & dark blue     \\ \hline
\end{tabular}
\end{center}
\end{table}

A non-solution with two sub-graphs is shown on the cover of the
instruction manual and is included below in Figure \ref{DemoNonSolutionLinedBlack}. A curve
separates the two halves of the puzzle. No pieces on either sub-graph
are connected to any pieces on the other sub-graph. (I have never
accidentally stumbled upon such a non-solution.)

The instruction manual claims that each posed problem has a unique
solution.  However, there used to be a published errata listing on the
official website \cite{arch} that indicated that the first version of
the toy's instruction booklet contained a few puzzles
could yield multiple solutions. This was fixed in future editions.
This does tell us, however, that it's possible to pose a solvable
puzzle with multiple solutions.  Naturally, if you are presented with
an entirely empty board, there are many solutions.

With those technical disclaimers out of the way, we can presume that the primary goal for a beginner puzzlist is to fit the pieces into the grid, and that a secondary check should be done at the end to ensure connectedness, just in case one should stumble across one of the disconnected non-solutions. 

\begin{figure}
\begin{center}
\includegraphics[width=3in]{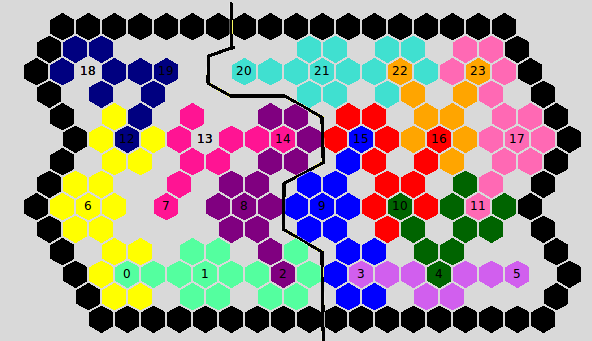}\\
\end{center}
\caption{A fitted placement which is not a solution.}
\label{DemoNonSolutionLinedBlack}
\end{figure}

\section{Constraints to piece placement }

There are 12 unique pieces. Of them, only one has flip symmetry. None have rotational symmetry. Each of the 12 pieces makes use of at least three ``spots'', either by using the inside of the spot or the outside, or both. (Complete circles effectively trap the inside, so this is modeled as occupation of the inner hexagon as well as the six outer hexagonal neighbors.)

In order to make a reasonable accounting for where the pieces lie, one ``ball'' or ``complete circle'' for each piece is designated as the anchoring point, although sockets would work also. (It is not important which ball or socket was selected. Of course, once an anchor is chosen, it's coded that way for all future calculations.) Each legal piece placement is therefore described uniquely by its flip (a binary choice available to all but one piece), its rotation angle (6 possibilities), and its anchor (one of 24 home anchor points).  There are, naturally, some illegal placements, because each ball or socket must rest on/in a ``spot'' on the board. In my model, I have an additional piece which is the boundary of the toy's play area. This way, I can test for collisions with the boundary to be sure the piece fits on the virtual board in the same way it would on the physical board. In my diagrams, this boundary piece is black.  Illegal placements are simply stopped in the algorithm by collision mechanisms just like other collisions between pieces.

The hexagonal diagram in Figure \ref{ReachablesAndHomes} shows the boundary, the 24 numbered ``spots'' (or ``anchors'') in yellow, and the unreachable locations in black. The toy's construction hides the fact that there are unreachable locations in the puzzle. There is no need to account for this separately in the code, as the anchoring process and allowable angles enforce this naturally.

As shown in Table \ref{puzzlepieces}, each piece has quite a few ways it can legally be placed on an empty board (but not all placements are necessarily strategic).  The three straight pieces should each give 144 legal placements. However, due to the mirror symmetry, the light purple piece has only 72 different placements. The four acute pieces have 180 legal placements, and the five obtuse pieces have 152.

\begin{figure}
\begin{center}
\includegraphics[width=3in]{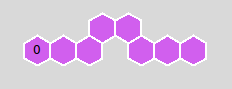}\\
\end{center}
\caption{The only piece with mirror symmetry.}
\label{mirror}
\end{figure}

Strategic concerns further reduce the options. For example, the light
purple piece is shown in Figure \ref{mirror} with its anchor point
labeled as ``0''.  There would be no reason to orient this piece as
shown if it were to anchor at spots 0 through 3. To do so would cause
its socket to point towards the boundary, rendering it unused.  It is
fine to leave a socket unused, however, if this piece could be placed
this way, it could be flipped, thus violating uniqueness.  There are
likewise a number of other useless moves which should never be used
from the other straight-family game pieces. We also know that sockets
of such useful placements must be in use when the placements are part
of the solution, in order to guarantee uniqueness.  These are
tabulated and filtered out in the code.

\begin{figure}
\begin{center}
\includegraphics[width=3in]{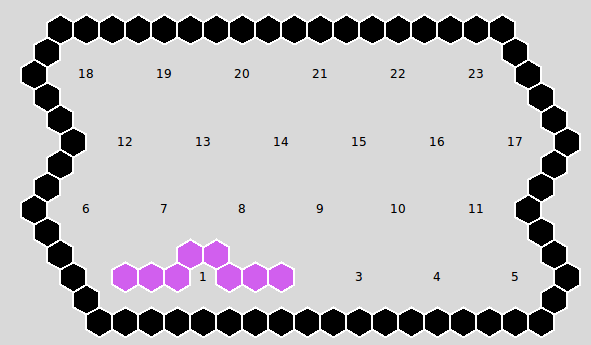}\\
\end{center}
\caption{A useless move.}
\label{uselessMove}
\end{figure}

\begin{figure}
\begin{center}
\includegraphics[width=3in]{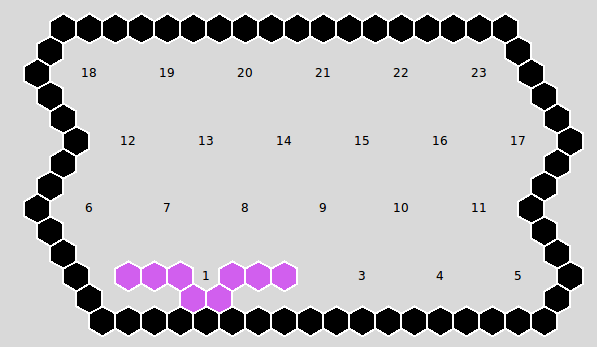}\\
\end{center}
\caption{A useful move.}
\label{usefulMove}
\end{figure}

\section{Orthographic coordinates }
To address a position in a hexagonal grid, orthographic coordinates are often used.  In a hexagonal grid, there are three families of parallel lines, thus we get ``tri-ordinates'' rather than ``coordinates''.  However, since the hex grid is in a plane, there is redundancy, and the third coordinate can be dropped computationally, although its use can lead to insights. 

\begin{figure}
\begin{center}
\includegraphics[width=6in]{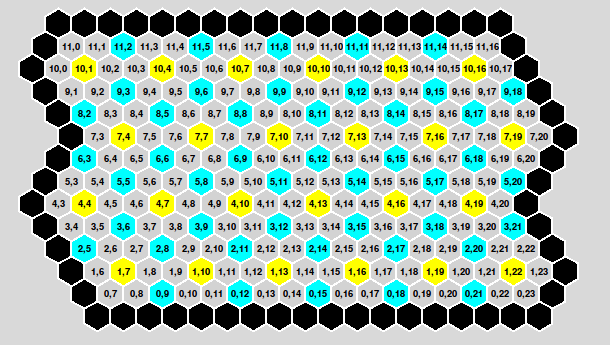}\\
\end{center}
\caption{Game board with  orthographic coordinates.}
\label{HexTemplateIQGame}
\end{figure}

Similar in nature to Figure \ref{ReachablesAndHomes}, Figure \ref{HexTemplateIQGame} shows the game board, but with orthographic coordinates. Horizontal lines have a common first coordinate, and parallel lines which slope upwards have a common second coordinate.  The third coordinate is not shown, but it is the sum of the other two.  Thus, elements on parallel lines which slope downwards have a common coordinate sum. The cyan elements are the non-reachable locations where were black in Figure \ref{ReachablesAndHomes}, and the yellow elements are the spots (or anchors) on the toy.

\begin{figure}
\begin{center}
\includegraphics[width=7in]{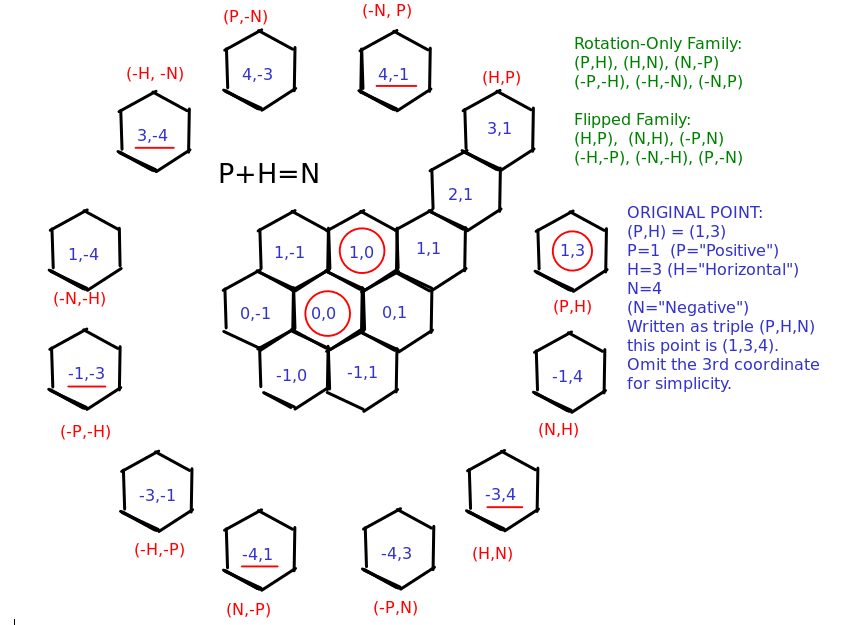}\\
\end{center}
\caption{How to flip and rotate in orthographic coordinates.}
\label{TriOrdinates}
\end{figure}

Consider Figure \ref{TriOrdinates}. The point (0,0,0) is labeled as (0,0) and indicates what might be the anchor point of a puzzle piece.  The point (1,3,4) is labeled as (1,3) and indicates a hex cell which is occupied by the puzzle piece in question. The  cell at (1,0) is circled to emphasize that it is upwards and to the right of (0,0).  The cell (1,3) is circled to emphasize that it is three units to the right of (1,0).  Thus, any cell can be reached by using just the two parallel line families P (for positive slope) and H (horizontal movements).  The third set of parallel lines is indicated by N (for negative slope).   Note again that since P+H=N, any two families may be used without loss of information.

Using (0,0) as the anchor point, the cell currently at (1,3) could end up at (4,-1,3), (3, -4, -1), (-1,-3,-4), (-4, 1,-3), or (-3, 4,1), that is, the underlined coordinates, depending on which legal rotation angle is selected.  If the piece is flipped, that cell could end up at (3,1,4), (4, -3, 1), (1, -4, -3), (-3,-1,-4), 
(-4,3,-1), or (-1,4,3).

Each puzzle piece is represented by a set of hex cell elements, each of which is translated by the same mechanism in order to determine the final cells of the game board which are occupied by it. There are 12 orientations for 11 pieces, and 6 orientations for the light purple piece. If there is a collision with the boundary, this is clearly an illegal placement.  Likewise, if a piece will collide with another piece, it may not be placed at the desired anchor with the desired orientation.

\section{``All Legal'' strategy }The first strategy one would probably consider is to simply try each available piece in every plausible spot.  Thus, the first thing tabulated by the algorithm is every placement that is legal for every piece, at the start of play.  For the first four ``Starter'' level puzzles, (in fact 7 of the first 11 puzzles) this alone yields the solution.  That is to say, as soon as you find that a piece fits into an open spot, it's in the correct final location for that piece.  This is where the puzzle starts to slowly ramp up in difficulty, as it builds confidence in the puzzlist.  However, it is also a great way to introduce kids and non-mathematician friends to a very fun puzzle.

The next thing a puzzlist might notice is that some (but certainly not all) pieces may have exactly one legal placement. Surely if this is the case for any piece, the solution must contain that placement.  The other pieces' options are now reduced.  This process iterates.  This basic process of elimination will solve 21 of the puzzles. Most of these are in the ``Starter'' level, but one is even in the ``Expert'' level!

\section{Pairwise cleanup }
To extend the process of elimination a bit, certain additional placements are eliminated. For each piece, there is a collection of possible legal placements upon starting the puzzle.   Consider one individual piece placement for a red piece, and another entire collection of placements for a green piece: If that one placement cannot be used with any of the placements for the second piece, clearly that placement cannot be used.  By this logic, quite a few placements are eliminated. This strategy alone solves all but one of the first 25 puzzles and a total of 29 puzzles overall.

\begin{figure}
\begin{center}
\includegraphics[width=3in]{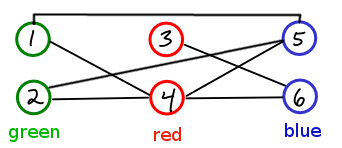}\\
\end{center}
\caption{Basic pairwise cleanup.}
\label{eliminatingOptions}
\end{figure}

As an example, let's take a puzzle which has six remaining legal placements, two for each of three pieces.  Vertices in Figure \ref{eliminatingOptions} indicate legal placements, and they are colored to indicate which piece they refer to. Edges represent compatible (non-intersecting) placement pairs.  Notice that placement 3 is not compatible with any of the placements for the green piece, so placement 3 can't be used.  Once placement 3 is removed from the graph, it is obvious that placement 4 must be in use, because it is the only remaining placement for the red piece. Similarly, placement 6 does not connect to any green placements, so it is discarded.

\section{The effects of basic cleanup techniques}
\label{basic}
\begin{figure}
\begin{center}
\includegraphics[width=3in]{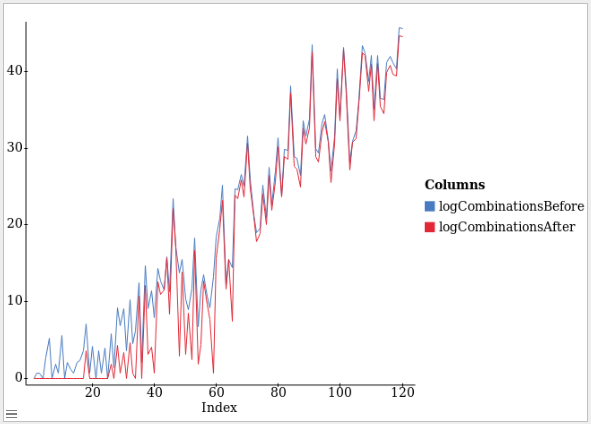}\\
\end{center}
\caption{Size of search space before and after basic cleanup (logscale).}
\label{LogCombinationsCleanup}
\end{figure}

After tabulating the available legal piece placements, basic pairwise cleanup and a few other easy tasks can dramatically reduce the search space.  These other basic tasks are: elimination the duplication from the symmetry of the light purple piece, and getting rid of the ``unwise'' placements. Figure \ref{LogCombinationsCleanup} compares before and after for the count of permitted piece placements, graphed on a log scale.

\begin{figure}
\begin{center}
\includegraphics[width=3in]{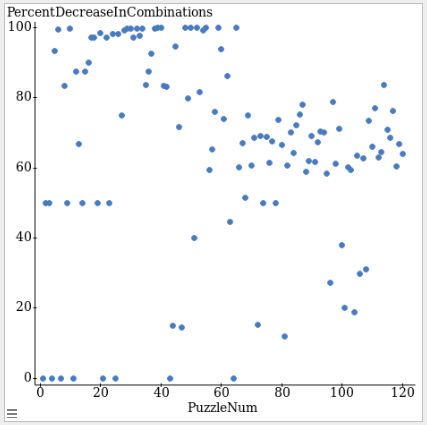}\\
\end{center}
\caption{Reduction of search space by percent.}
\label{CleanupPercentImprove}
\end{figure}
Even expressed on a log scale, the reduction in possible moves does not look very impressive, and it's even less exciting without the log scale.  However, when expressed as total combinations of piece placements (Figure \ref{CleanupPercentImprove}) rather than just the number of piece placements options, and when the results are expressed as a percent rather than raw values, it's clear that there is quite a bit of value in doing even this little bit of cleanup.  For the easiest problems, it seems to help quite a bit, but less so for puzzles after about number 60.  

\section{Heatsort}
\label{heatsort}
Since a naive way to brute force a solution is to simply try all possible permutations of legal placements for each piece, it would make sense to try to prioritize the list somehow.  One way this can be done is to superimpose all the possible legal placements onto a single grid and see which locations on the base grid are harder to reach.  Although it's true that not all spots must be covered, it is reasonable, for example, to prioritize trying ``the one piece that can go there'' before trying to leave a position empty.  That being said, this is something of an Achilles' Heel for this algorithm.

\begin{figure}
\begin{center}
\includegraphics[width=3in]{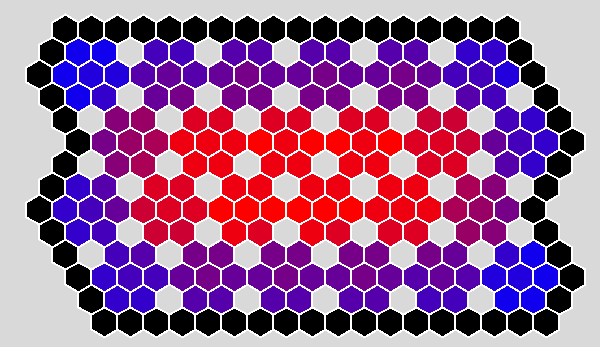}\\
\end{center}
\caption{The edges are not as easily occupied as the center of the board.}
\label{colorHeat_starting}
\end{figure}

\begin{figure}
\begin{center}
\includegraphics[width=3in]{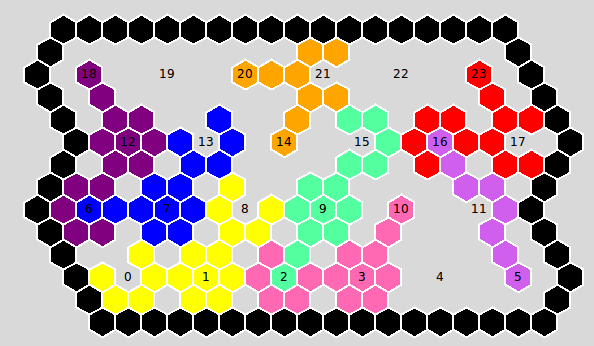} \includegraphics[width=3in]{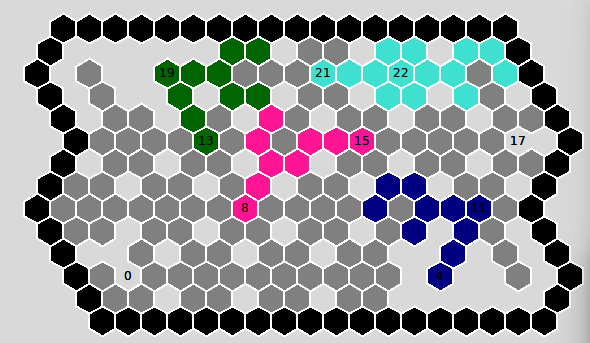} \includegraphics[width=3in]{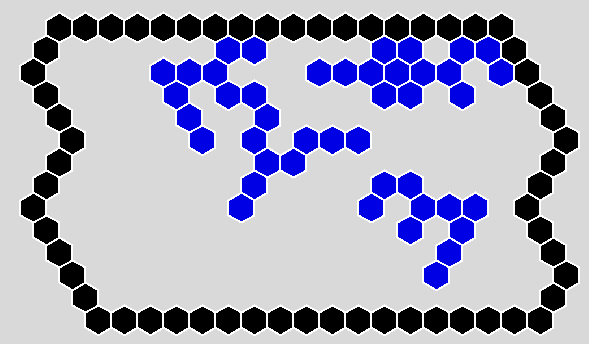}\\
\end{center}
\caption{Puzzle 25, its solution, and its heatmap.}
\label{puz25}
\end{figure}

\begin{figure}
\begin{center}
\includegraphics[width=3in]{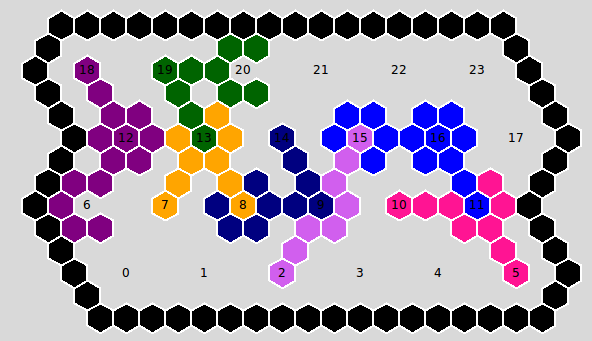} \includegraphics[width=3in]{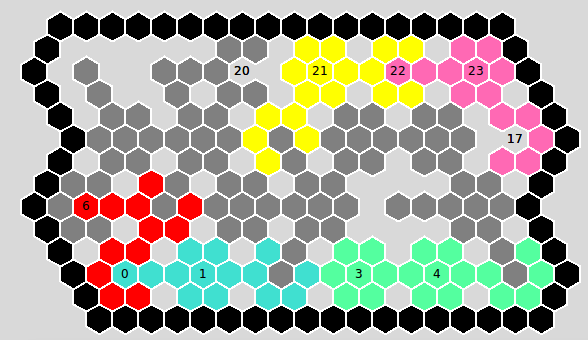} \includegraphics[width=3in]{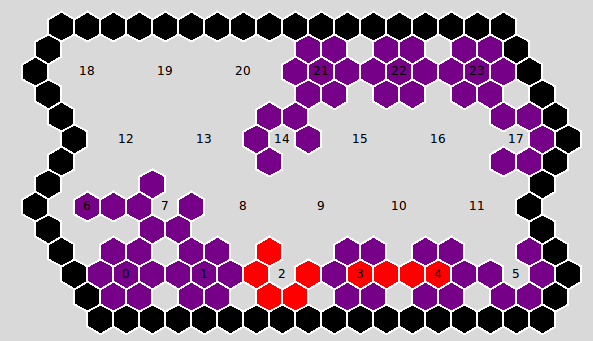}\\
\end{center}
\caption{Puzzle 29, its solution, and its heatmap.}
\label{puz29}
\end{figure}

\begin{figure}
\begin{center}
\includegraphics[width=3in]{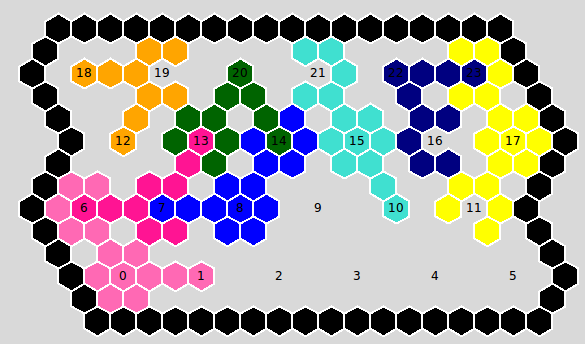}  \includegraphics[width=3in]{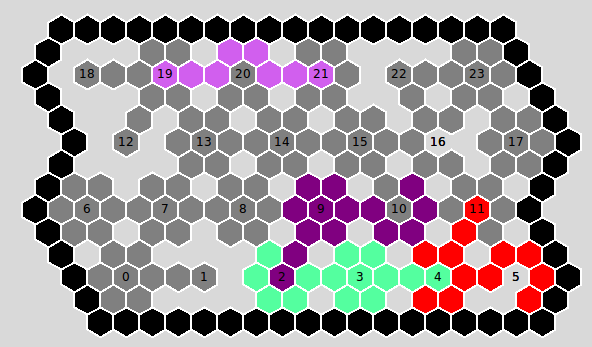}  \includegraphics[width=3in]{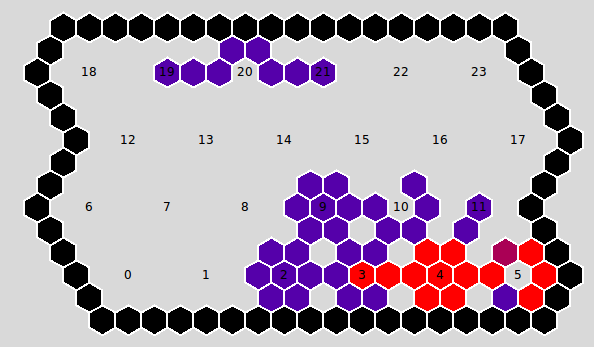}\\
\end{center}
\caption{Puzzle 40, its solution, and its heatmap.}
\label{puz40}
\end{figure}

\begin{figure}
\begin{center}
\includegraphics[width=3in]{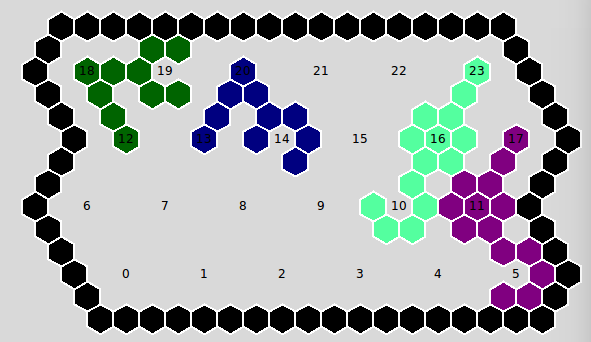} \includegraphics[width=3in]{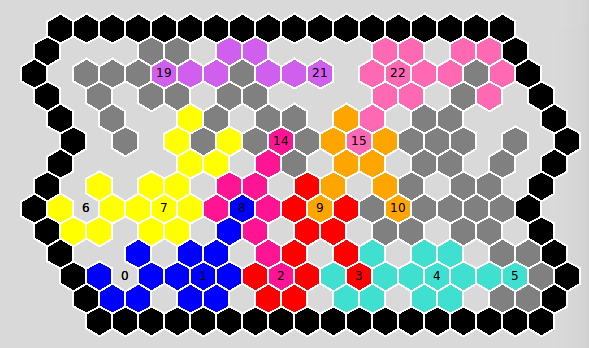} \includegraphics[width=3in]{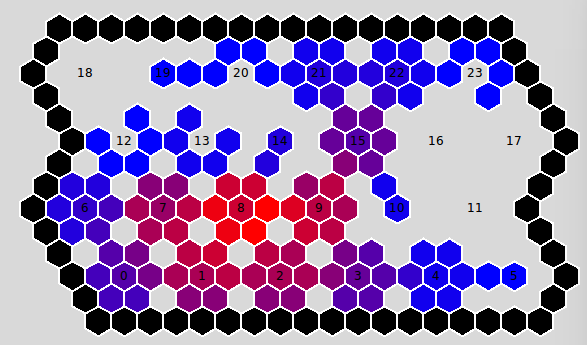}\\
\end{center}
\caption{Puzzle 85, its solution, and its heatmap.}
\label{puz85}
\end{figure}

The heatmap in Figure \ref{colorHeat_starting} reflects the superposition of all possible legal moves, with red indicating a more reachable place and blue indicating a harder to reach place. For four puzzles (Figures \ref{puz25}, \ref{puz29}, \ref{puz40}, and \ref{puz85}), the puzzle, the solution with the puzzle grayed out, and the heatmap are shown.

Heatmap 25 is all one color, because every element of this heatmap has the value 1.  Each piece has exactly one legal placement, and hence the puzzle is trivial.

The heatmap for Puzzle 29 has mostly purple elements, but the bits of red indicate elements which can be covered by more than one puzzle piece.

In Puzzle 40, the heatmap is still giving very clear directives for us to place a couple of pieces with priority.  For example, the light purple piece must be placed at the top of the puzzle. However, for Puzzle 85, you see more color variations, and the heatmap looks more like the overall heatmap for the entire puzzle.  Thus, the algorithm will first try something such as placing the correct light purple, light pink, and aqua pieces down and then continue to look for
the solution to the remaining puzzle with 6 missing pieces.

\begin{figure}
\begin{center}
\includegraphics[width=6in]{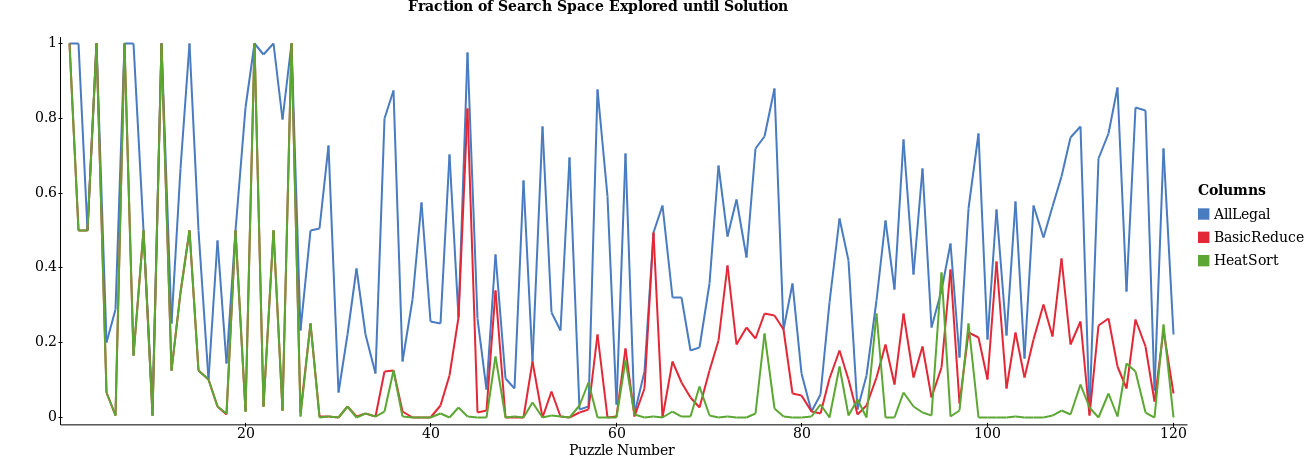}\\
\end{center}
\caption{Fraction of search space explored, reduced by basic techniques, then by heatsort.}
\label{FractionSearchSpaceExplored}
\end{figure}

The benefits of adding heatmap sorting to our basic introductory phase of the sorting can be shown in Figure \ref{FractionSearchSpaceExplored}. Since a solution method (Section \ref{linkage}) was eventually successful in finding all 120 puzzle solutions, it was possible to compute how many puzzle combinations would have had to have been tested from a given sorted order of legal piece placements, in order to reach the solution.  The number of combinations (as shown) grows with increasing polynomial complexity, so this graphic can be misleading.  Still, it is interesting to see what fraction of the combinations are washed out by the various techniques.  

For the simpler puzzles, not only is there less benefit to doing any sorting, there simply aren't many combinations to begin with.  For example, in Puzzle 6, the 3rd of 15 combinations to be tested will be found to work.  After basic reductions, only one viable combination remains. For Puzzle 19, the puzzle would be initially solved after 178 of 1232 combinations. That drops to the 12th combination after basic reductions and heat sort. Puzzle 119, however, starts out with around  $7\times10^{19}$ combinations and has to search through about $5\times10^{19}$ of them before reaching the solution.  After cleaning up and heat sort, this drops to $2\times10^{19}$.  What a bargain!  Puzzle 120 appears to be fairly ``easy'' with only $6\times10^{19}$ combinations, dropping to a mere $2\times10^{16}$. So the fraction of search space to be explored does not say much about computability.

\section{Search space restriction by linkage requirement }
\label{linkage}
Even with the reductions discussed, and even after using a heatsort,
the search space is still enormous.  But the addition of linkage requirements to the algorithm ended up making the difference in computability.  If not for
this requirement being added, there might still be computers running! Recall that the
number of pieces to be placed keeps growing throughout the puzzle
levels, although not monotonically.  To solve a puzzle with $n$ missing
pieces is therefore $O(n!)$ (which is consistent with \cite{hearn}
p.13) for a naive search mechanism, even after using clever sorting
tricks. Let me pause to make a disclaimer about the order of this
puzzle and then retract it.  Since the maximal number of missing
pieces is 12, one can argue that this puzzle is never higher than
$O(n^{12} )$ and cannot rightfully claim its place in the elite
community of factorial problems.  However, while this is true, it is
clear to me that this puzzle was capped at precisely that point at
which the factorial growth starts to kick in computationally, making
it indeed a part of the discussion.

To reduce this horrific computational explosion, a wiser strategy (and more compatible with what a human puzzlist would do) would be to restrict each successive piece placement to a piece which would actually link to one of the existing pieces already in the puzzle.  This makes sense because a puzzle can be completed by only selecting from the available, currently-linkable pieces. 

Please recall the connected graph structure requirement from the puzzle's instructions in Section \ref{mechanics}. The solution has to be a connected graph, so the solution itself can be thought of as its own spanning tree. Loops are certainly permitted but never required, and they do not impede the analysis in any way.

\begin{figure}
\begin{center}
\includegraphics[width=3in]{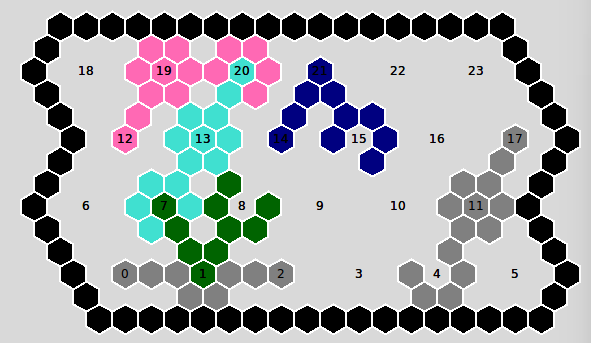}\\
\end{center}
\caption{Two grayed out solution pieces for Puzzle 92, one of which is currently linkable.}
\label{Explain92}
\end{figure}

As a specific example, consider the puzzle in Figure \ref{Explain92} (previously also shown in Figure \ref{pubpuz92}).  The light purple piece can be anchored immediately at spot 2 as shown in gray in Figure \ref{Explain92}, and it's part of the solution.  Another placement, which is also part of the solution and is also grayed out, cannot yet be anchored to its home at spot 17, because it is not currently linkable. 

\begin{figure}
\begin{center}
\includegraphics[width=3in]{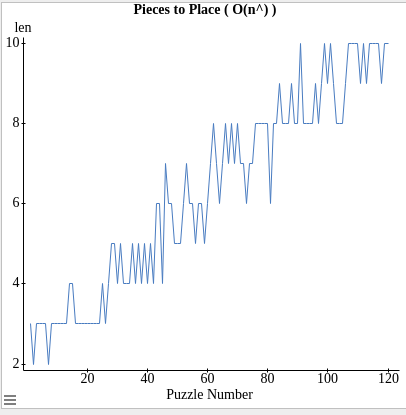}\\
\end{center}
\caption{Pieces to place, graphed by puzzle number.}
\label{PiecesToPlace_Onxx}
\end{figure}

\begin{figure}
\begin{center}
\includegraphics[width=3in]{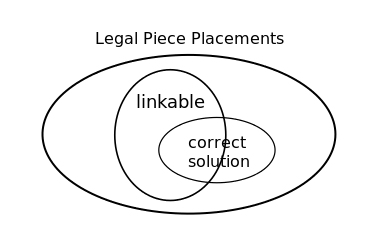}\\
\end{center}
\caption{Relationship between all legal placements, correct placements, and currently linkable pieces.}
\label{LegalPieceVenn}
\end{figure}

\begin{figure}
\begin{center}
\includegraphics[width=3in]{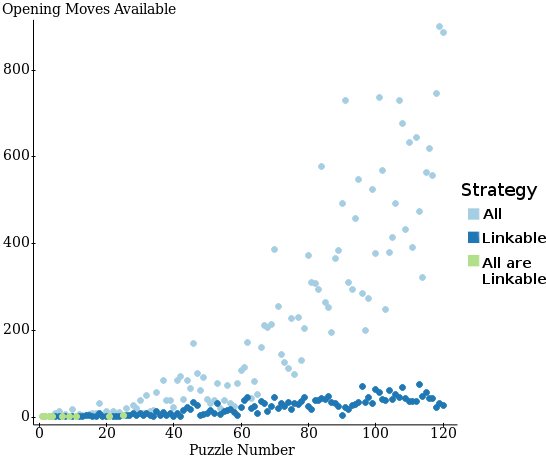}\\
\end{center}
\caption{Total number of legal and linkable moves at the start of each challenge.}
\label{F24}
\end{figure}

While this reduces the search space, it makes for a more complex process.  Consider the first piece to be placed in any puzzle, especially a harder one. There will be many available placements, and all placements in the final solution are immediately available for you to select.  However, if you force each successive placement to be a linked placement, this dramatically reduces the number of pieces you may  place at any given moment, however, if you are already on the right track, you are promised at least one fruitful placement is available to you to make.  This makes an apples-to-oranges comparison difficult, but it has to be done nonetheless.  Figure \ref{F24} indicates how many placements are available at the start of each puzzle.  Note that for harder puzzles, this can be over 800.

Figure \ref{F24} also shows the number of original placements which are also linkable to the original puzzle.  The highest values are around 80.  For the easy puzzles, all available placements are linked.

This focus on only linkable placements actually reduces the problem dramatically, to the point where it is finally computable.  If linkage conditions are not used, the problem quickly explodes.  This effort is not without cost, and for easier puzzles, it is not worth the time computationally.  The transition, however, is rather brutal. (Computed on FUJITSU T734 laptop.)

To give an idea of the pruning that happens at runtime, consider Puzzle 29. As the algorithm searches for possible moves, at one point it finds 23 legal placements, only 5 of which satisfy the linkability criterion.  At some other point, it finds 20 legal placements, 4 of which satisfy it. This is shown in Table \ref{reduction}. To give an idea of what happens to a random sample of puzzles of various difficulties, this table gives a partial listing of some of the decision reductions made at various steps by the algorithm.

For the puzzles selected here, the reductions range roughly between $80\%$ and $99\%$, and the solution times were quite fast, with Puzzle 80 taking about 24 seconds and Puzzle 90 taking about 2 minutes.

Figure \ref{F24} shows the number of first moves available, per puzzle, in light blue, and it shows the number of linkable first moves available in dark blue.  For the easier problems, all moves are linkable and are shown in light green.  However, for the most difficult problems, there is a huge differential, going from numbers close to 900 down to numbers closer to 80.  

\begin{table}[]
    \begin{center}
  \caption{Comparision of runtimes for two methods for a few select puzzles running on a FUJITSU T734.}
  \label{transition}
  \begin{tabular}{|l|l|l|l|}
\hline
\textbf{puzzle} & \textbf{missing pieces} & \textbf{search without linkability} & \textbf{search with linkability} \\ \hline
47              & 7                                 & (gave up!)                                     & 2.5 seconds                                 \\ \hline
48              & 6                                 & 0.001 seconds                                  & 1.5 seconds                                 \\ \hline
49              & 5                                 & 2 seconds                                      & 1.5 seconds                                 \\ \hline
50              & 5                                 & 0.001 seconds                                  & 0.7 seconds                                 \\ \hline
51              & 5                                 & 0.02 seconds                                   & 0.8 seconds                                 \\ \hline
52              & 6                                 & 0.001 seconds                                  & 1.3 seconds                                 \\ \hline
53              & 7                                 & 170 seconds                                    & 2.7 seconds                                 \\ \hline
  \end{tabular}
  \end{center}
\end{table}

\begin{table}[]
     \begin{center}
  \caption{Reduction ratios for linkability criteria for selected puzzles.}
  \label{reduction}
\begin{tabular}{|l|l|l|}
\hline
\multicolumn{1}{|c|}{puzzle} & \multicolumn{1}{c|}{number of missing pieces} & \multicolumn{1}{c|}{ratios on search branches} \\ \hline
29                           & 5                                             & 5/23  4/20                                     \\ \hline
38                           & 4                                             & 4/41                                           \\ \hline
60                           & 6                                             & 22/129  19/95  11/47                           \\ \hline
80                           & 8                                             & 25/407  16/324  13/267  7/98                   \\ \hline
90                           & 8                                             & 5/550  5/442  11/359  7/79  3/33  5/85.....    \\ \hline
\end{tabular}
\end{center}
\end{table}

\section{Transition between the two methods }
The algorithm finally used to solve the puzzles was basically
recursive with various guiding strategies, including heatsort (Section
\ref{heatsort}) and use of linkability (Section \ref{linkage}), but when
it gets down to either 3 or 4 remaining pieces, it searches
exhaustively.  The transition to exhaustive search can be set to any
level, but experience shows that 3 or 4 remaining pieces is the best
transition point, and it doesn't matter too much which of those is
selected.  Recall Table \ref{transition} where it is demonstrated that
this transition between the two basic strategies is rather abrupt. The
reader may wonder why 5 or 6 is not used as the transition level,
because of the values in Table \ref{transition}, but in the table, we
were only looking for known solutions, which can be found much more
quickly than exhastively searching an entire space of solutions to
verify that none exist in any branch.  For example, we were very lucky
in Puzzle 53 to reach the solution to a Puzzle with 7 missing pieces
in only 170 seconds.  We were able to do this only because the first
guess was correct due to heatsort and thus the solution was found
quickly. This was not the case, for example, for Puzzle 47.

\begin{figure}
\begin{center}
\includegraphics[width=7in]{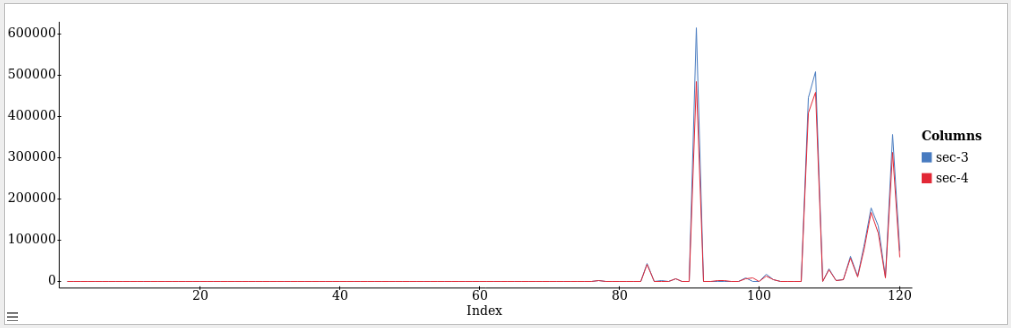}\\
\end{center}
\caption{Computation time per challenge, in seconds.}
\label{computeSec}
\end{figure}

\begin{figure}
\begin{center}
\includegraphics[width=7in]{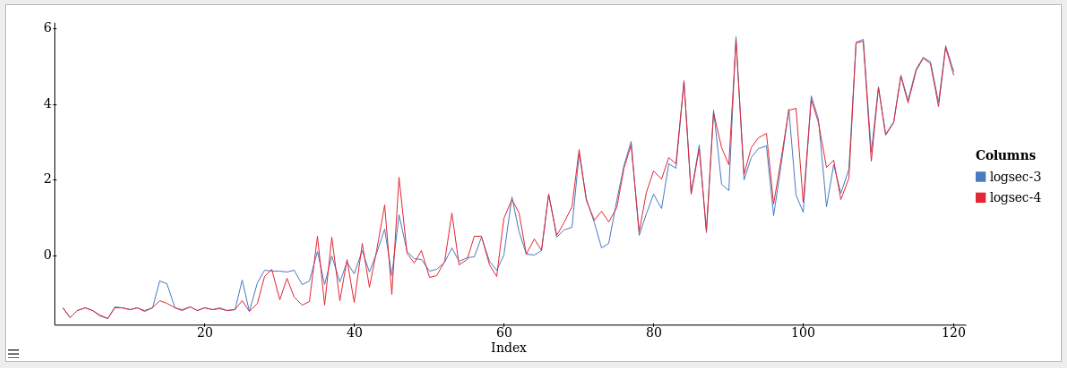}\\
\end{center}
\caption{Computation time per challenge, with logscale.}
\label{computeLogSec}
\end{figure}

The two graphs in Figures \ref{computeSec} and \ref{computeLogSec}  give computation time in seconds, and the log-seconds for transitioning to exhaustive search at either 3 or 4 remaining pieces in a puzzle. 
All puzzles were ultimately solved. (There can only be a solution or a hung computer. We are guaranteed to eventually find the solution.)  The longest puzzle took about 5 1/2 to 7 days to solve on an Intel core i7 - 3820 8-core processor running at 3.6 GHz with 32 GB of RAM, depending on which transitioning strategy was used.  Notice that the spikes which look really curious Figure \ref{computeSec} are hardly even noticeable in Figure \ref{computeLogSec}, due to the factorial nature of the growth of move options in the game.

\section{Counting balls and sockets }
One thing that I found myself tracking by hand was the amount of physical space available in the puzzle.  Spots may be left empty according to the instructions.
But how many can we leave empty? Recall that the board contains 24 spots. 

Of the 12 pieces, there are these broad categories: ``2 Balls and a Socket'',  ``1 Ball and 2 Sockets'',  ``2 Sockets and a full Circle'', and ``1 Ball, 1 full Circle, 1 Socket''.  If you could break up those pieces, you would get 16 Balls, 14 Sockets, and 6 Circles.  Circles take up the anchoring spot as well as the 6 hexagonal spots which surround the anchor.  Most balls can mate with most sockets, but there are a couple of exceptions which I can safely ignore for this analysis.  (The dark green and dark blue acute pieces each have an ``elbow'' style ball.  These can only fit into the light purple socket.)

Therefore, we can use the 24 spots in just one of the following three ways:
2 unmatched balls, 14 ball/socket sets, and 2 empty spots
3 unmatched balls, 1 unmatched socket, 13 ball/socket sets, and 1 empty spot
4 unmatched balls, 2 unmatched sockets, 12 ball/socket sets, and no empty spots.

A more robust algorithm can look for permanent isolations of this
sort, and once the maximum number of unmatched balls or empty sockets
is reached, no more should be permitted, causing a return to the
previous recursion level.  Balls and sockets were determined to be
permanently unmatched when the heat map indicated that no remaining
pieces could be placed with them.

It's harder to make a general statement about this part of the
algorithm, because I only did a few side-by-side timed comparisions.
In a similar pattern to what we saw in Section \ref{basic}, the harder
the puzzle, the more time is saved by doing this analysis.  For easier
problems, the time sunk to do the calculations is essentially
wasted. To give a few examples, Puzzle 75 took 53 seconds on the
FUJITSU T734, down from 66 seconds without testing for unreachable balls
and sockets.  Puzzle 80 took 24 seconds either way. However, Puzzle 82
went from 813 seconds without the test down to 492 seconds with
them. For all these cases, the transition to brute force was made when
there were 3 pieces remaining.

\section{Other insights}Many of the smaller details of the algorithm have been left to the reader's imagination, but there are a few insights worthy of sharing. Once you calculate the set of legal piece placements for a puzzle (which is actually pretty fast), the original puzzle as well as the boundary may be disposed of.  Since all collisions with the original puzzle and the boundary are not longer possible, the only collisions must be between solution pieces. This cuts down the computational complexity.  This insight has logical implications, such as ``well let's just tabulate all possible collections of 3 or 4 or 5 pieces!''  So far, this instinct and many other presumably mathematically correct insights have led to many failed solution computation attempts, because of computational complexities that are not immediately obvious.  There is always hope that there is yet another mathematical insight that can be made that will further reduce the computational complexity of these searches.

\section{Insights not incorporated in the algorithm }One possible improvement to the algorithm would be to attempt to implement what seems obvious to a human puzzlist but might be difficult for a computer.  Sometimes, the puzzle appears to have an obvious split, or is highly likely to have a split, which causes the puzzle to appear as two (or perhaps more) separate puzzles sharing the same pool of pieces.  For a human, this adds information, but in terms of coding, I have really no idea how to put that part of my intuition into the algorithm. 

There are certainly many more improvements which could be made to the algorithm, but the primary objective was to generate advice for the would-be addict of this puzzle.  Towards that end, I feel I've gathered enough evidence to support my conclusion.

\section{Conclusion }Even after having ``solved'' all 120 puzzles using a combination of search-space reductions and brute force, I still find myself in awe of the puzzle.  For me, I can conclude that it is a wholesome toy that encourages logical thinking, because the puzzlists must continue to update their strategies as they advance through the stages in order to complete the challenges.  A human cannot solve the puzzle simply by trying all possible combinations, except for at the very earliest stages.  Since the puzzle is inexpensive and aesthetically pleasing, I encourage would-be puzzlists to go ahead and give it a try!

\end{document}